\documentclass[12pt]{article}
\usepackage{amsfonts, amssymb, amsmath,color}

\usepackage{enumerate} 

\let\Horig\H

\usepackage{scalerel}

\usepackage[colorlinks=true, urlcolor=blue,linkcolor=blue, citecolor=blue]{hyperref}

\usepackage{graphics,graphicx,amsmath}

\newcommand{\R}{\mathbb{R}}

\newcommand{\E}{\mathbb{E}}

\newcommand{\N}{\mathbb{N}}

\newtheorem{prop}{Proposition}[section]
\newtheorem{lemma}[prop]{Lemma}

\newtheorem{corollary}[prop]{Corollary}
\newtheorem{theorem}[prop]{Theorem}

\def\({\left(}
\def\){\right)}

\def\[{\left[}
\def\]{\right]}
\def\real{{\mathord{\mathbb R}}}
\def\inte{{\mathord{\mathbb N}}}
\def\Dom{\mathrm{Dom}}
\def\Var{\mathrm{Var}}
\newcommand{\p}{\mathbb{P}}

\newenvironment{Proof}{\removelastskip\par\medskip
\noindent{\em Proof.} \rm}{\penalty-20\null\hfill$\square$\par\medbreak}

\allowdisplaybreaks

\numberwithin{equation}{section}

\begin{document}
\title{
\huge
Normal approximation for generalized $U$-statistics and weighted random graphs 
} 

\author{Nicolas Privault\thanks{Division of Mathematical Sciences,
Nanyang Technological University, SPMS-MAS-05-43, 21 Nanyang Link
Singapore 637371. e-mail: {\tt nprivault@ntu.edu.sg}. 
}
\and 
Grzegorz Serafin\thanks{Faculty of Pure and Applied Mathematics, Wroc{\l}aw University of Science and Technology, Ul. Wybrze\.ze Wyspia\'nskiego 27, Wroc{\l}aw, Poland.
e-mail: {\tt grzegorz.serafin@pwr.edu.pl}.}}

\maketitle

\vspace{-0.4cm}

\begin{abstract} 
 We derive normal approximation bounds in the Wasserstein distance
 for sums of weighted $U$-statistics, based on a general distance bound for
 functionals of independent random variables of arbitrary distributions. 
 Those bounds are applied to normal approximation for
 the combined weights of subgraphs in the Erd{\H o}s-R\'enyi
 random graph, extending the graph counting results of \cite{BKR}
 to the setting of graph weighting. 
 Our approach relies on a general stochastic analytic
 framework for functionals of independent random sequences. 
\end{abstract}
 
\baselineskip0.7cm

\section{Introduction}
The Malliavin calculus has been applied to the derivation 
of approximation bounds by the Stein and Chen-Stein methods
on the Wiener space \cite{nourdinpeccati}
and on the Poisson space \cite{utzet2}, \cite{thale},
see also \cite{nourdin3}, \cite{privaulttorrisi4}, \cite{reichenbachs}, \cite{reichenbachsAoP}, \cite{krokowskicosa} for the case of discrete Bernoulli sequences.
Recently, a different Malliavin framework for Stein approximation
has been introduced in \cite{PS}, with application to normal
approximation in the Wasserstein distance for weighted $U$-statistics
of the form 
$$ 
\sum_{ \substack{k_1,\ldots ,k_n\in\N_0\\k_i\neq k_j \text{ if }i\neq j}}
 b_{k_1} \cdots b_{k_n} Z_{k_1} \cdots Z_{k_n} 
$$
 where $\N_0=\{0,1,2,\ldots \}$,
 $(Z_k)_{k\geq 1}$ is an i.i.d. sequence of random
 variables, 
 and $(b_k)_{k\geq }$ is a sequence of real coefficients,
 based on stochastic analysis for functionals of a countable number
of uniformly distributed random variables, see 
\cite{prebub}. 
This completes the bounds for the Kolmogorov distance
obtained in e.g. Theorem~3.1 of \cite{CS} for non weighted $U$-statistics,
see also \cite{goetze} in the quadratic case. 
\\

Our goal in the present paper is two-fold.
First, we extend in Theorem~\ref{thm:Var<Sum} the Stein approximation bounds of 
\cite{PS} from {multiple stochastic integrals to finite sums of multiple stochastic integrals, which can be viewed as polynomial functionals
in independent random variables with
arbitrary distributions, or as generalized weighted $U$-statistics, see Proposition \ref{prop:sum=I}.}
Furthermore, in Proposition~\ref{thm:dWI-1} we obtain a general
Wasserstein distance bound for
functionals of independent random variables 
as a consequence of Proposition~3.3 in \cite{PS}.
\\

Second, we show that those results can be applied to
the central limit theorem for the
convergence of renormalized weight counts in large random graphs.
For this, we consider the Erd{\H o}s-R\'enyi random
graph $\mathbb{G}(n,p)$, 
introduced by Gilbert \cite{G} in 1959 and popularized
in \cite{ER},
which is constructed by independently retaining any edge in 
the complete graph $K_n$ on $n$ vertices with probability $p\in(0,1)$.
Denote by 
$N^G_n$ the random variable counting number of  subgraphs  (not necessarily induced ones) of $\mathbb{G}(n,p_n)$
that are isomorphic to a fixed graph $G$.
Necessary and sufficient conditions for the
asymptotic normality of the renormalization
$$\widetilde{N}^G_n:=\frac{N^G_n-\E[N^G_n]}{\sqrt{\Var[N^G_n]}}.
$$
 of $N^G_n$ have been obtained in \cite{rucinski}, where it is shown that  
\begin{align} \label{eq:conv}
  \widetilde{N}^G_n\stackrel{\mathcal D}{\longrightarrow}\mathcal{N}\ \mbox{ iff }
  \ \ np_n^\beta\rightarrow \infty \ \mbox{ and }\ n^2(1-p_n)\rightarrow\infty,
\end{align}
as $n$ tends to infinity, where
$\mathcal{N}$ represents the standard normal distribution, 
 $\beta=\beta(G) :=\max \{e_H/v_H \ : \ H\subset G\}$ 
and $e_H$, $v_H$ respectively denote
the numbers of edges and vertices in the graph $H$.
Such results have been improved via
explicit convergence rates obtained in \cite{BKR}
as
 \begin{equation} 
  \label{ddww}
  d_W \big(\widetilde{N}^G_n,\mathcal{N}\big)\leq C\((1-p_n)\min_{\substack{ H\subset G\\e_H\geq1}} n^{v_H}p_n^{e_H} \)^{-1/2},
\end{equation} 
 where $\mathcal{N}$ represents the standard normal distribution, $C>0$
 is a constant depending on $G$, and $d_W$ is the Wasserstein distance
$$
 d_W (X,Y):
 =\sup_{h\in\mathrm{Lip}(1)}|\mathrm{E}[h(X)]-\mathrm{E}[h(Y)]|,
$$
 between the laws of random variables $X$, $Y$,
 where $\mathrm{Lip}(1)$ denotes the class of real-valued
 Lipschitz functions with
 Lipschitz constant less than or equal to $1$.
Kolmogorov distance bounds have also been obtained for triangle counting, 
see \S~3.2.1 of \cite{nross}, and \cite{reichenbachsAoP}, 
using the Malliavin approach to 
the Stein method for discrete Bernoulli sequences. 
Those rates have been improved in \cite{roellin2}, and extensions to the counting of arbitrary subgraphs that yield the bound \eqref{ddww} for the Kolmogorov distance have recently been obtained in \cite{PS2}, based on distance bounds for sums of discrete multiple integrals and weighted $U$-statistics, as well as in \cite{xiaofang}, \cite{zhuosongzhang}. 
\\ 

Here, our stochastic analytic framework allows us to assign 
an independent sample of a random nonnegative weight $X$ to every edge
in $\mathbb{G}(n,p_n)$,
and to consider the combined weights of subgraphs
instead of counting them. Precisely, we define a weight of a graph as a sum of weights of its edges.  Next, by $W^G_n$ we denote the combined weight of
subgraphs in $\mathbb{G}(n,p_n)$ that are isomorphic to
a fixed graph $G$ and its renormalization
\begin{equation}
  \label{wg} 
\widetilde{W}^G_n:=\frac{W^G_n-\E[W^G_n]}{\sqrt{\Var[W^G_n]}}.
\end{equation} 
In Theorem~\ref{thm:main}
we show, as an application of Corollary~\ref{thm:Var<Sum},
 that when $G$ is a graph without isolated vertices, 
 we have
\begin{equation} 
  \label{eq:main}
  d_W \big(\widetilde{W}^G_n,\mathcal{N}\big) \leq C
  \frac{\sqrt{\E[(X-\E[X])^4]}+(1-p_n)(\E[X])^2}{\Var [X]+(1-p_n)(\E[X])^2}
  \((1-p_n)\min_{\substack{ H\subset G\\e_H\geq1}}
  n^{v_H}p_n^{e_H} \)^{-1/2},
\end{equation} 
where $C >0$ is a constant depending only on $e_G$,
which recovers \eqref{ddww} in the case of a
deterministic weight given by $X:=1/e_G$.
When $X$ is a fixed random variable this
also yields the sufficient condition
\begin{align*}
 \(np_n^\beta\rightarrow \infty \ \mbox{ and }\ n^2(1-p_n)\rightarrow\infty\) \Longrightarrow \widetilde{N}^G_n\stackrel{\mathcal D}{\longrightarrow}\mathcal{N},
\end{align*}
for the convergence of $\widetilde{W}^G_n$ to
 the standard normal distribution (cf. \eqref{eq:conv}), 
which follows from the equivalence 
$$\(np_n^\beta\rightarrow \infty \ \mbox{ and }\ n^2(1-p_n)\rightarrow\infty\)\Longleftrightarrow (1-p_n)\min_{\substack{ H\subset G\\e_H\geq1}}
  n^{v_H}p_n^{e_H}\rightarrow\infty. $$
To derive the bound \eqref{eq:main} we
apply Proposition~\ref{thm:dWI-1}
to combined subgraph weights
$W^G_n$ represented as
finite sums of multiple stochastic integrals, 
see Lemma~\ref{lem:W=sumI}.  
Our results are then specialized to
a class of graphs satisfying a certain balance condition,
which includes triangles, complete graphs and trees
as particular cases.
\\

We note that other types of random functionals on graphs, 
such as graph weights defined as products of edge weights,
or the number of vertices of a given degree,
admit representations as sums of multiple integrals
or weighted $U$-statistics, and can be treated by this approach. 
\\

\noindent 
 This paper is organized as follows. 
 In Section~\ref{s1} we recall 
 the framework of \cite{prebub} for the construction
 of random functionals of uniform random variables,
 together with the construction of derivation operators. 
 In Section~\ref{s2} we derive normal Stein approximation
 bounds for general functionals and for sums of
 multiple stochastic integrals. 
 In Section~\ref{s3} we show that combined graphs weights
 can be represented as sums of multiple stochastic
 integrals,
 and
 derive distance bounds
 for the renormalized
 weights of graphs in $\mathbb{G}(n,p_n)$ that are isomorphic to
 a fixed graph $G$. The Appendix Section~\ref{sec-app}
 contains some technical results exploited in the paper.  
\section{Functionals of uniform random sequences}
\label{s1}
\subsubsection*{Stochastic integrals}
Given $(U_k)_{k\in \inte}$
an i.i.d. sequence
of $[-1,1]$-valued uniform random
variables on a probability space
$(\Omega , {\cal F}, P)
=
( [-1,1]^\inte , {\cal F}, P)$
let the jump process $(Y_t)_{t\in \real_+}$ be defined as 
$$Y_t : =\sum_{k=0}^\infty\mathbf1_{[2k+1+U_k,\infty)}(t),
  \qquad t\in \real_+. 
$$
  Denoting by $({\cal F})_{t\in \real_+}$
  the filtration generated by $(Y_t)_{t\in \real_+}$, and letting  
  $$\tilde{\cal F}_t: ={\cal F}_{2k},\qquad 2k\leq t<2k+2,
  \qquad k \in \inte, 
  $$
  the compensated stochastic integral
  $$
  \int_0^{\infty}u_t d(Y_t-t/2)
  $$
  with respect to the compensated point process 
  $(Y_t - t/2)_{t\in \real_+}$ can be defined for square-integrable
  $(\tilde{\cal F}_t)_{t\in \real_+}$-adapted processes
  $(u_t)_{t\in \real_+}$ by the isometry relation 
  \begin{equation}
    \label{eq:uv}
    \E\[\int_0^{\infty} \hskip-0.1cm
    u_t d(Y_t-t/2)
    \int_0^{\infty}
    \hskip-0.1cm
    v_t d(Y_t-t/2)\]
=
\E\[\int_0^{\infty} \hskip-0.1cm 
u_t \(v_t - \frac{1}{2} \sum_{k=0}^\infty\mathbf1_{(2k,2k+2]}(t)\int_{2k}^{2k+2}
\hskip-0.1cm
v_r dr
\)\frac{dt}{2}\],
\end{equation} 
  see \cite{prebub}, 
  where $(u_t)_{t\in \real_+}$ and $(v_t)_{t\in \real_+}$ are
square-integrable $(\tilde{\cal F}_t)_{t\in \real_+}$-adapted processes.
This also implies the bound
\begin{equation}
\nonumber 
  \E\[\left(
\int_0^{\infty} \hskip-0.1cm
    u_t d(Y_t-t/2)
\right)^2 \]
\leq 
\frac{1}{2}
\E\[\int_0^{\infty} \hskip-0.1cm 
| u_t |^2 dt\],
\end{equation}
where $(u_t)_{t\in \real_+}$ is square-integrable and 
$(\tilde{\cal F}_t)_{t\in \real_+}$-adapted.
\subsubsection*{Multiple stochastic integrals}
Let $\widehat{L}^p(\R_+^n)$
denote the space of symmetric functions that
are $pth$ integrable on $\R_+^n$,
$p\geq 1$, and vanish outside of 
$$
\Delta_n : = \bigcup_{
  0\leq k_i \not= k_j \atop
  1\leq i \not= j \leq n
}
      [2k_1,2k_1+2]\times \cdots \times [2k_n,2k_n+2], 
$$
      equipped with the norm
      $$
      \Vert f_n \Vert_{\widehat{L}^p(\R_+^n)} :
      =
      \Vert f_n \Vert_{L^p(\R_+^n ,(dx/2)^{\otimes n})}
      =
      \frac{1}{2^{n/p}}
      \Vert f_n \Vert_{L^p (\R_+^n , (dx)^{\otimes n})}
      ,
      \qquad f_n \in \widehat{L}^p(\R_+^n).
      $$
      Given $f_n\in \widehat{L}^1(\R_+^n) \cap \widehat{L}^2(\R_+^n)$,
      $n\geq 1$, we define the multiple stochastic integral $I_n(f_n)$
      as 
\begin{eqnarray}
  \label{page5} 
  \lefteqn{ I_n(f_n)
    := \sum_{r=0}^n
    \frac{(-1)^{n-r}}{2^{n-r}} {{n}\choose{r}} 
  }
  \\
  \nonumber
  & & 
  \sum_{ \substack{k_1,\ldots ,k_r\in\N_0\\k_i\neq k_j \text{ if }i\neq j}}
  \int_0^\infty\cdots \int_0^\infty f_n (2k_1+1+U_{k_1}
  ,\ldots ,2k_r+1+U_{k_r} ,y_1,\ldots ,y_{n-r})dy_1\cdots dy_{n-r}
  \\
\nonumber 
  &= & n!\int_0^\infty \int_0^{t_n}\cdots \int_0^{t_2}f_n(t_1,\ldots ,t_n)d(Y_{t_1}-t_1/2)\cdots d(Y_{t_n}-t_n/2), 
\end{eqnarray} 
 where $\N_0=\{0,1,2,\ldots \}$, with in particular 
$$
I_1(f_1)
:= \sum_{k=0}^\infty f_1 (2k+1+U_k )
- \frac{1}{2} \int_0^\infty f_1 (t) dt 
 = \int_0^\infty f_1(t) d(Y_t-t/2), 
$$
 for $f_1\in L^1(\R_+)\cap L^2(\R_+)$.
 The multiple stochastic integrals
 $(I_n(f_n))_{n\geq 1}$ form a 
family of mutually orthogonal centered random variables 
with the bound 
\begin{equation}  \label{eq:I^2<}
  \E \big[( I_n(f_n))^2\big]\leq n!\left\|f_n\right\|^2_{\widehat{L}^2(\R_+^n,dx/2)},
  \qquad n\geq 1, 
\end{equation} 
cf. \eqref{eq:uv} above and Propositions~4 and 6 of \cite{prebub},
which allows one to extend the definition of $I_n(f_n)$ to all
$f_n \in \widehat{L}^2(\R_+^n)$.
If in addition we have 
\begin{equation}
\label{ass:int0}
\int_{2k}^{2k+2}f_n (t,*)dt=0, \qquad k \in \inte, 
\end{equation} 
then $I_n(f_n)$ satisfies 
the isometry and orthogonality relation 
\begin{equation}
\label{eq:EI^2}
\E\[ I_n(f_n) I_m(f_m) \]=
           {\bf 1}_{\{ n = m \} }
           n!\langle f_n , f_m \rangle_{\widehat{L}^2(\R_+^n,dx/2)},
           \quad f_n \in \widehat{L}^2(\R_+^n), \quad
           f_m \in \widehat{L}^2(\R_+^m), 
\end{equation}
see \cite{prebub}, page~589.
Moreover, every $F \in L^2(\Omega)$ admits the chaos decomposition 
\begin{equation}
  \label{eq:decom}
F=\E[F] + \sum_{n=1}^\infty I_n(f_n),
\end{equation}
for some sequence $(f_n)_{n\geq 1}$ of functions
in $\widehat{L}^2 (\R_+^n)$, $n\geq 1$, see Proposition~7 of \cite{prebub}. 
Note that under the condition \eqref{ass:int0}
the sequence $(f_n)_{n\geq 1}$ is unique
in $\widehat{L}^2 (\R_+^n)$ due to the isometry
relation~\eqref{eq:EI^2}.

\subsubsection*{Finite difference operator} 
Consider the finite difference operator $\nabla$
defined on multiple stochastic integrals $F = I_n(f_n)$ as
\begin{equation}
  \nonumber 
  \nabla_tF := F \circ \Phi_t - \frac{1}{2}
  \int_{2\lfloor t/2\rfloor}^{2\lfloor t/2\rfloor+2}F\circ \Phi_s ds,
  \qquad t \in \real_+, 
\end{equation}
where $\Phi_t : \Omega \longrightarrow \Omega$ 
is defined by 
$$\Phi_t(\omega) : =\(U_1(\omega ),\ldots ,U_{\lfloor t/2\rfloor-1}(\omega ),t-2\lfloor t/2\rfloor-1,U_{\lfloor t/2\rfloor+1}(\omega ),\ldots \),
\quad \omega \in \Omega, \quad t\in \real_+. 
$$
 The operator $\nabla$ admits an adjoint operator $\nabla^*$ given by 
$$
 \nabla^* \left( I_n(g_{n+1}) \right) : = I_{n+1}({\bf 1}_{\Delta_{n+1}}
 \tilde{g}_{n+1}), 
$$ 
where $\tilde{g}_{n+1}$ is the symmetrization of
$g_{n+1} \in \widehat{L}^2(\R_+^n)\otimes L^2(\R_+)$
in $n+1$ variables.
The operator $\nabla$ is closable with domain
$$
\Dom ( \nabla ) =
\big\{ F\in L^2(\Omega ) : \E [ \Vert \nabla F \Vert^2_{L^2(\real_+)} ]
< \infty \big\}
  , 
  $$
and we have the duality relation (integration by parts)
\begin{equation}
 \label{dr} 
\E[\langle \nabla F , u \rangle_{\widehat{L}^2(\R_+)} ] 
=
\E[ F \nabla^*(u)], \qquad F \in \Dom ( \nabla ),
\end{equation}
 for $u$ in the domain $\Dom ( \nabla^* )$ of $\nabla^*$,
 see Proposition~8 of \cite{prebub}. 
 Although the operator $\nabla$ does not satisfy the chain rule of derivation,
 it can be easily applied to multiple stochastic integrals,
 as for any $f_n\in \widehat{L}^2(\R_+^n)$ we have 
\begin{align}\label{eq:nablaext}
\nabla_t I_n(f_n)=n I_{n-1}(f_n(t,*))-n\int_{2\lfloor t/2\rfloor}^{2\lfloor t/2\rfloor+2} I_{n-1}(f_n(s,*))ds, \qquad t\in \real_+, 
\end{align}
see Proposition~2.1 in \cite{PS}.
In particular, under the condition \eqref{ass:int0} we have the equality
$$
 \nabla_t I_n(f_n) = n I_{n-1}\(f_n(t,*)\), \qquad
 t\in \real_+, 
$$
 see Proposition~10 of \cite{prebub}.
 The Ornstein-Uhlenbeck operator $L : =-\nabla^*\nabla$ satisfies 
\begin{equation*}
  L I_n(f_n) =-\nabla^*\nabla I_n(f_n) =-n I_n(f_n),
  \qquad
  f_n\in \widehat{L}^2(\R_+^n), 
\end{equation*}
where $f_n$ satisfies \eqref{ass:int0}.
By \eqref{eq:decom}
the operator $L$ is well defined, invertible on centered random variables 
$F\in L^2(\Omega)$, and its inverse operator $L^{-1}$ is given by
\begin{equation}
  \nonumber 
L^{-1} I_n(f_n)  =-\frac1n I_n(f_n),\qquad n\geq1,
\end{equation}
where, due to Proposition~\ref{prop:Ifoverline} below,
$f_n$ does not have to satisfy \eqref{ass:int0}.
Note that $(-L)$ is a positive operator and its square root  $(-L)^{-1/2}$ takes the form
\begin{align*}(-L)^{1/2} I_n(f_n)  =\sqrt n I_n(f_n),\qquad n\geq1.\end{align*}
\subsubsection*{Stein approximation bound} 
The next result is a consequence of
 Proposition~3.3 in \cite{PS}. 
 \begin{prop}
   \label{thm:dWI-1}
   Let $X\in \Dom ( \nabla )$ be such that $\E[X]=0$. 
   We have 
\begin{eqnarray} 
\label{dwf-2} 
d_W(X , \mathcal{N}) & \leq & 
\left|1 - \E [ X^2 ]     \right|        +
    \sqrt{\Var\[ \langle\nabla_\cdot X, -\nabla_\cdot L^{-1}X\rangle_{\widehat{L}^2(\R_+)}\]}
    \\
    \nonumber
    & & + \,2 \sqrt{\E\[ | (-L)^{-1/2} X |^2\]
   \int_0^\infty  \E \left[ 
| \nabla_t X |^4 
 \right]\,\frac{dt}2}\\
 \label{thm:dWI-2}
 &\leq&|1-\E[X^2]|+
\sqrt{\Var\[ \langle\nabla_\cdot X, -\nabla_\cdot L^{-1}X\rangle_{L^2(\R_+)}\]}
  \\
    \nonumber
    & &+\,2 \sqrt{\E[X^2]
   \int_0^\infty  \E \left[ 
| \nabla_t X |^4 
\right]\,\frac{dt}2}
.
\end{eqnarray} 
\end{prop}
\begin{Proof}
The inequality \eqref{thm:dWI-2} follows from \eqref{dwf-2} by Proposition  \ref{prop:sm} with $F=L^{-1}X$, so it is enough to prove \eqref{dwf-2}. Proposition~3.3 in \cite{PS} states that 
\begin{align} 
\nonumber 
  d_W(X,\mathcal{N}) \leq\,& \E\[\left|1-\frac{1}{2} \langle\nabla_\cdot X, -\nabla_\cdot L^{-1}X\rangle\right|\]
\\
\nonumber
& 
 + \frac{1}{2}
    \E \left[ 
\int_0^\infty 
| \nabla_t L^{-1} X | 
| \nabla_t X |^2 dt
 \right] +  
\frac{1}{4}
    \E \left[ 
\int_0^\infty 
| \nabla_t L^{-1} X | 
\int_{2\lfloor t/2\rfloor}^{2\lfloor t/2\rfloor+2}
 | \nabla_s X |^2 ds dt
 \right].
\end{align}
We will estimate each of the three terms on the left-hand side. First,  taking $F=X$ and $u=\nabla L^{-1}X$ in \eqref{dr}, we get
\begin{align} 
\nonumber 
&\E\[\left|1-\frac{1}{2} \langle\nabla_\cdot X, -\nabla_\cdot L^{-1}X\rangle\right|\]
    \\
\nonumber
    & \leq 
    \E\[\left|1
    -\frac{1}{2} \E \left[ \langle\nabla_\cdot X, -\nabla_\cdot L^{-1}X\rangle
      \right]
    \right| 
    \]
    +
    \E\[\left|\frac{1}{2} \langle\nabla_\cdot X, -\nabla_\cdot L^{-1}X\rangle
    -\frac{1}{2} \E \left[ \langle\nabla_\cdot X, -\nabla_\cdot L^{-1}X\rangle
      \right]
    \right|\]
    \\
\nonumber
    & \leq 
   \left|1 - \E [ X^2 ]     \right|        +
    \sqrt{\Var\[ \langle\nabla_\cdot X, -\nabla_\cdot L^{-1}X\rangle_{\widehat{L}^2(\R_+)}\]}
.
\end{align} 
Next, for $F=L^{-1}X$ in \eqref{eq:sm}, we obtain
\begin{align*}
  \E \left[ \int_0^\infty | \nabla_t L^{-1} X |^2 dt\right] 
   &=  2 \E\[ | (-L)^{-1/2} X |^2\] .
\end{align*}
Consequently, the Cauchy-Schwarz inequality gives us
\begin{align*}
\E \left[ \int_0^\infty | \nabla_t L^{-1} X | | \nabla_t X |^2 dt\right] &\leq \sqrt{\E \left[ \int_0^\infty | \nabla_t L^{-1} X | ^2 dt\right]\E \left[ \int_0^\infty | \nabla_t X |^4 dt\right]}\\
&\leq 2
\sqrt{\E\[ | (-L)^{-1/2} X |^2\]}
\sqrt{\E \left[ \int_0^\infty | \nabla_t X |^4\, \frac{dt}2\right]},
\end{align*}
and
\begin{eqnarray*}
  \lefteqn{
    \! \! \! \! \! \! \! \! \! \! \! \! \! \! \! \! \!  \! \! \!  \! \! \! 
    \E \left[ \int_0^\infty | \nabla_t L^{-1} X | 
      \int_{2\lfloor t/2\rfloor}^{2\lfloor t/2\rfloor+2}| \nabla_s X |^2 ds dt \right]
=\sum_{k=0}^\infty  \E \left[ \int_{2k}^{2k+2} | \nabla_t L^{-1} X | dt
  \int_{2k}^{2k+2}| \nabla_s X |^2  ds \right]
  }
  \\
&=& \sqrt{ \E \left[\sum_{k=0}^\infty \(\int_{2k}^{2k+2} | \nabla_t L^{-1} X | dt\)^2\]
\E\[\sum_{k=0}^\infty\(\int_{2k}^{2k+2}| \nabla_s X |^2\)^2  ds \right]}\\
  &\leq & 2
  \sqrt{ \E \left[ \int_0^\infty | \nabla_t L^{-1} X |^2 dt\]\E\[\sum_{k=0}^\infty\int_{2k}^{2k+2}| \nabla_s X |^4  ds \right]}\\
   &\leq &4 \sqrt{\E\[ | (-L)^{-1/2} X |^2\]}
   \sqrt{ \E\[\int_0^\infty| \nabla_s X |^4 \, \frac{dt}2 \right]},
\end{eqnarray*}
and we conclude  \eqref{dwf-2}. 
\end{Proof}

\section{Normal approximation for weighted $U$-statistics} 
\label{s2}
In this section we consider
generalized weighted $U$-statistics of order $n \geq 1$
of the form
$$ 
 \sum_{\substack{k_1,\ldots ,k_n\in\N_0\\k_i\neq k_j \text{ if }i\neq j}}
 f_n (2k_1+1+U_{k_1},\ldots ,2k_n+1+U_{k_n} ),
 $$
 where $\N_0=\{0,1,2,\ldots \}$. The next proposition gives 
the multiple stochastic integral expansion of
such extended weighted $U$-statistics. 
\begin{prop}\label{prop:sum=I}
 Given $f_n \in \widehat{L}^2(\R_+^n)$ we have
\begin{align} 
\nonumber 
 \sum_{ \substack{k_1,\ldots ,k_n\in\N_0\\k_i\neq k_j \text{ if }i\neq j}}
  f_n (2k_1+1+U_{k_1},\ldots ,2k_n+1+U_{k_n} ) =  
   \sum_{r=0}^n  I_k\(f_n^{(r)}\), 
\end{align} 
 where 
$$f_n^{(r)}=(x_1,\ldots ,x_k)=
\frac{1}{2^{n-r}}
     {{n}\choose{r}} 
     \int_{\R^{n-r}_+}f_n(x_1,\ldots ,x_r,y_1,\ldots ,y_{n-r})dy_1\cdots dy_{n-r},
          $$
$r=0,1,\ldots , n$.
\end{prop}
\begin{Proof}
Formula \eqref{page5} gives us for $1\leq m\leq n$
\begin{align} 
\nonumber 
 a_m:&=(-2)^mI_m\(\int_{\R^{n-m}}f_n\(*,y_1,\ldots,y_{n-m}\)dy_1\cdots dy_{n-m}\)  = \sum_{r=0}^m
    (-1)^r {{m}\choose{r}} b_r,
\end{align}
where 
\begin{align*}
  b_r=2^r \sum_{\substack{k_1,\ldots ,k_r\in\N_0\\k_i\neq k_j \text{ if }i\neq j}}
  \int_0^\infty\cdots \int_0^\infty f_n (2k_1+1+U_{k_1}
  ,\ldots ,2k_r+1+U_{k_r} ,y_1,\ldots ,y_{n-r})dy_1\cdots dy_{n-r}.
\end{align*} 
Hence, by binomial inversion, we have $b_m=\sum_{r=1}^m(-1)^r{m \choose r}a_r$, $1\leq m\leq n$. In particular, 
\begin{align} 
\nonumber 
&\sum_{ \substack{k_1,\ldots ,k_n\in\N_0\\k_i\neq k_j \text{ if }i\neq j}}
  f_n (2k_1+1+U_{k_1},\ldots ,2k_n+1+U_{k_n} )
\\\nonumber
&=2^{-n}b_n=2^{-n}\sum_{r=1}^n(-1)^r{n \choose r}a_r\\\nonumber
&\ \  = 
  \sum_{r=0}^n
  {{n}\choose{r}} 
  \frac{1}{2^{n-r}}
  I_r\(\int_{\R^{n-r}}f_n\(*,y_1,\ldots,y_{n-r}\)dy_1\cdots dy_{n-r}\)
  \\
  \nonumber
  &\ \ =  \sum_{r=0}^nI_r\(f_n^{(r)}\),
\end{align} 
as required.
\end{Proof}
 In particular, under the condition \eqref{ass:int0}
 the multiple stochastic integral $I_n(f_n)$
coincides with the weighted $U$-statistic of order $n$ and we have  
\begin{align}
 I_n(f_n)&=\sum_{ \substack{k_1,\ldots ,k_n\in\N_0\\k_i\neq k_j \text{ if }i\neq j} }f_n (2k_1+1+U_1,\ldots ,2k_n+1+U_n). 
\end{align}
In the next corollary we obtain
a Wasserstein distance bound
for sums of multiple stochastic integrals
by combining Propositions~\ref{prop:fg<f^2+g^2}
and \ref{thm:dWI-1} with the multiplication formula \eqref{eq:mf}. First, let us introduce the following $\star$-notation: for $0\leq l\leq k\leq n\wedge m$
we define the contraction 
$f_n\star_{k}^lg_m$ of
$f_n \in \widehat{L}^2(\R_+^n)$ and
$g_m \in \widehat{L}^2(\R_+^m)$ as 
\begin{align}\label{starnotation}
 & f_n\star_{k}^lg_m(x_1,\ldots , x_{k-l},y_1,\ldots,y_{n-k},z_1,\ldots , z_{m-k}) 
  \\\nonumber
   & := \frac{1}{2^l} 
\int_{\R^l_+} f_n(w_1,\ldots,w_l, x_1,\ldots , x_{k-l},y_1,\ldots,y_{n-k}) \\\nonumber
&\ \ \ \ \ \ \ \ \ \ \ \ \ \ \times g_m(w_1,\ldots , w_l,x_1,\ldots , x_{k-l},z_1,\ldots , z_{m-k})dw_1\cdots dw_l,
\end{align}
           and we let
           $f_n \hskip0.1cm \widetilde{\star}_{k}^l g_m$
           denote the symmetrization
           \begin{eqnarray*}
             \lefteqn{ 
               \! \! \! \! \! \! \! \! \! 
               f_n \hskip0.1cm \widetilde{\star}_{k}^l \hskip0.1cm g_m(x_1,\ldots , x_{n+m-k-l})
             }
             \\
             & := &
             \frac{\mathbf1_{\Delta_{m+n-k-l}}(x_1,\ldots , x_{n+m-k-l})
              }{(m+n-k-l)!}\sum_{\sigma\in S_{m+n-k-l}}f_n\star_{k}^lg_m(x_{\sigma(1)},\ldots ,x_{\sigma(m+n-k-l)}),
\end{eqnarray*} 
where $S_n$, $n\geq1$, denotes the set of all permutations of the set $\{1,\ldots ,n\}$.
  \begin{theorem}
  \label{thm:Var<Sum}
  For any $X\in L^2 (\Omega )$ written as a sum $X=\sum_{k=1}^n I_k(f_k)$
  of multiple stochastic integrals where  
  $f_k\in \widehat{L}^2(\R_+^k)$ satisfies \eqref{ass:int0},
  $k=1,\ldots ,n$,
  we have 
\begin{eqnarray*}
  \lefteqn{
    d_W(X , \mathcal{N})  \leq   
 \left|1 - \E [ X^2 ]     \right|     
  }
  \\
 & &
 +
  C_n \sqrt{
    \sum_{0\leq l< i\leq n}\left\| f_{i}\star_{i}^l f_{i}\right\|^2_{
      L^2(\R_+^{i-l})}+\sum_{1\leq l< i\leq n}
    \(\left\| f_{i}\star_{l}^l f_{i}\right\|^2_{
      L^2(\R_+^{2(i-l)})}+\left\| f_{l}\star_{l}^l f_{i}\right\|^2_{
      L^2(\R_+^{i-l})}\)}, 
\end{eqnarray*}
for some $C_n>0$. 
\end{theorem}
\begin{Proof}
 Given that 
\begin{align*}
  \nabla_tX = \sum_{k=0}^{n-1} (k+1) I_k \(f_{k+1}(t,\cdot)\), \quad
  \mbox{and}
  \quad
  \nabla_tL^{-1}X = \sum_{k=0}^{n-1} I_k \(f_{k+1}(t,\cdot)\), 
  \end{align*}
 the multiplication formula \eqref{eq:mf} shows that 
 \begin{equation}
   \label{XX}
   (\nabla_t X)^2=\sum_{0\leq i\leq j < n }\sum_{k=0}^{i}\sum_{l=0}^kc_{i,j,k,l}I_{i+j-k-l}
   \big(
   f_{i+1}(t,\cdot)
   \hskip0.1cm \widetilde{\star}_{k}^l f_{j+1}(t,\cdot)\big)
   \end{equation} 
   and
   \begin{equation}
\label{XLX}
\nabla_t X\nabla_t L^{-1}X=\sum_{0\leq i\leq j < n }\sum_{k=0}^{i}\sum_{l=0}^kd_{i,j,k,l}I_{i+j-k-l}
\big( f_{i+1}(t,\cdot)
\hskip0.1cm \widetilde{\star}_{k}^l
f_{j+1}(t,\cdot)\big),
\end{equation} 
   for some $c_{i,j,k,l}, d_{i,j,k,l}\geq0$.
   Next, by \eqref{eq:I^2<} and \eqref{XX} we get
\begin{align}\nonumber
  \int_0^\infty  \E \left[ | \nabla_t X |^4  \right]\,\frac{dt}2&\leq
  C
  \sum_{0\leq i\leq j < n }\sum_{k=0}^{i}\sum_{l=0}^k\int_0^\infty
  \big\|f_{i+1}(t,\cdot)\widetilde\star_{k}^lf_{j+1}(t,\cdot)\big\|^2_{\widehat{L}^2(\R_+^{i+j-k-l})}dt
  \\
  \nonumber
&\leq C
  \sum_{0\leq i\leq j < n }\sum_{k=0}^{i}\sum_{l=0}^k \left\|f_{i+1} \hskip0.1cm {\star}_{k+1}^lf_{j+1} \right\|^2_{{L}^2(\R_+^{i+j-k-l+1})}\\\label{aux1}
&\leq C
  \sum_{1\leq i\leq j\leq n}\sum_{k=1}^{i}\sum_{l=0}^{k-1} \left\|f_{i} \star_{k}^lf_{j} \right\|^2_{{L}^2(\R_+^{i+j-k-l})}, 
\end{align}
where $C>0$ is a constant depending on $n$.
Furthermore, from \eqref{XLX} it follows that 
\begin{eqnarray*}
  \lefteqn{
    \langle\nabla_\cdot X, -\nabla_\cdot L^{-1}X\rangle-\E\[\langle\nabla_\cdot X, -\nabla_\cdot L^{-1}X\rangle\]
  }
  \\
  &=&
  \frac{1}{2}
  \int_0^\infty\sum_{0\leq i\leq j < n }\sum_{k=0}^{i}\sum_{l=0}^kd_{i,j,l,k}\mathbf1_{\{i=j=k=l\}^c}I_{i+j-k-l}\big(
  f_{i+1}(t,\cdot)\widetilde \star_{k}^lf_{j+1}(t,\cdot)\big) dt,
\end{eqnarray*}
thus we get 
\begin{align}\nonumber
\Var&\[\langle\nabla_\cdot X, -\nabla_\cdot L^{-1}X\rangle\]\\\nonumber
&\leq C'\sum_{0\leq i\leq j < n }\sum_{k=0}^{i}\sum_{l=0}^k\mathbf1_{\{i=j=k=l\}^c}\left\|\int_0^\infty f_{i+1}(t,\cdot)\star_{k}^lf_{j+1}(t,\cdot)dt\right\|^2_{{L}^2(\R_+^{(i+j-k-l)})}\\\nonumber
&= C''
\sum_{0\leq i\leq j < n }\sum_{k=0}^{i}\sum_{l=0}^k\mathbf1_{\{i=j=k=l\}^c}\left\| f_{i+1}\star_{k+1}^{l+1} f_{j+1}\right\|^2_{{L}^2(\R_+^{i+j-k-l})}\\\label{aux2}
&= C''
\sum_{1\leq i\leq j\leq n}\sum_{k=1}^{i}\sum_{l=1}^k\mathbf1_{\{i=j=k=l\}^c}\left\| f_{i}\star_{k}^l f_{j}\right\|^2_{{L}^2(\R_+^{i+j-k-l})},
\end{align}
for some constants $C',C''>0$ depending only $n$. Applying \eqref{aux1} and \eqref{aux2} to \eqref{thm:dWI-2}, we get
\begin{equation*}
  d_W(X , \mathcal{N} )\leq
  \left|1 - \E [ X^2 ] \right|  +
  C''' \sqrt{
    \sum_{1\leq i\leq j\leq n}\sum_{k=1}^{i}\sum_{l=0}^k\mathbf1_{\{i=j=k=l\}^c}\left\| f_{i}\star_{k}^l f_{j}\right\|^2_{{L}^2(\R_+^{i+j-k-l})}},
\end{equation*}
for some $C'''>0$ depending on $n$. Next, by the inequality \eqref{eq:l<k}, all the components where $0\leq l<k\leq i,j$, are dominated by those where $0\leq l<k=i=j$, and also, by the inequality \eqref{eq:l=k},  the ones where $1\leq k=l<\min\{i,j-1\}$, are dominated by the components where $1\leq l=k<i=j$.
Finally, the components for $1\leq k=l=i<j$ remain unchanged. 
\end{Proof}

\section{Application to weighted random graphs} 
\label{s3}
\noindent
In this section we present an application of results from the previous section to the Erd{\H o}s-R\'enyi random
graph $\mathbb{G}(n,p)$ and
to the renormalization $\widetilde{W}^G_n$ 
of the combined weight $W_n^G$
of subgraphs of the random graph that are isomorphic to
a fixed graph $G$, see \eqref{wg}.

{
In order to simplify the notation we write $a_n\lesssim b_n$
for two sequences $a_n$ and $b_n$
whenever there exist a constant $C$ depending only on $G$ such that $a_n<Cb_n$ for all $n\in\N$. Furthermore, if $a_n\lesssim b_n$ and $b_n\lesssim a_n$
then we write $a_n\approx b_n$. Finally, by writing
$H\sim K$ we mean that the two graphs $H$ and $K$ are isomorphic.}
In Proposition \ref{prop:VarW} we provide estimates of the variance of $W_n^G$, which is crucial when dealing with the renormalization. 
\begin{prop}\label{prop:VarW}
  The variance of $W^G_n$ admits the asymptotic form 
  \begin{align}
    \label{eq:VarW}
    \Var\big[W^G_n\big]\approx
    \big(
    \Var[X]+(1-p_n)(\E[X])^2
    \big)
    \max_{\substack{H\subset G \\ e_H\geq1}} n^{2v_G-v_H}p_n^{2e_G-e_H}. 
\end{align}
\end{prop}
\begin{Proof}
 We follow the lines of the proof of Lemma~3.5 in \cite{JLR}
 by extending the argument to nonnegative random weights distributed as $X$.
 We note that
$$
W^G_n = \sum_{G'\sim G} S_{G'}, 
$$ 
 where the sum is over all graphs $G' \subset K_n$ 
 which are isomorphic to $G$, and
 $S_{G'}$ is the sum of the weights of edges in $G'$
 if $G'$ belongs to $\mathbb{G}(n,p_n)$,
 and zero otherwise, i.e. denoting by $X_1,\ldots ,X_{e_G}$ the random weights of edges of $G'$, we have
 $$S_{G'}:=\bold{1}_{\{G'\in \mathbb{G}(n,p_n)\}}\sum_{i=1}^{e_G}X_i.$$
  Then, we get
\begin{align*}
\Var\big[W^G_n\big]&=\sum_{G', G'' \sim G}\text{Cov}(S_{G'},S_{G''})\\
&=\sum_{G', G''  \sim G \atop
  {\rm with\ a\ common\ edge}}
\big(\E[S_{G'}S_{G''}]-\E[S_{G'}]\E[S_{G''}]\big)
\\
&\approx
\sum_{\substack{H\subset G\\e_H\geq1}}\sum_{G'\cap G'' \sim H
  \atop
  G', G'' \sim G
}\big(\E[S_{G'}S_{G''}]-\E[S_{G'}]\E[S_{G''}]\big).
\end{align*}
    {
      For a fixed $G'\sim G$ we clearly have
      $$
      \E[S_{G'}]=\p(G'\in \mathbb{G}(n,p_n))\sum_{i=1}^{e_G}\E[X_i]=e_Gp^{e_G}\E[X].
      $$
      In order to calculate $\E[S_{G'}S_{G''}]$ for $G',G''\sim G$ and $G'\cap G''\sim H$,
      let us denote by $X_1, \ldots ,X_{e_H}$ the
      weights of edges of $G'\cap G''$ and by $X'_1,\ldots ,X'_{e_G-e_H}$ and $X''_1,\ldots ,X''_{e_G-e_H}$ weights of edges of $G'\backslash G''$ and $ G''\backslash G'$, respectively. Then, we have 
\begin{align*}
 &\E[S_{G'}S_{G''}]=\p\(G',G''\in  \mathbb{G}(n,p_n)\)\E\[\(\sum_{i=1}^{e_H}X_i+\sum_{i=1}^{e_G-e_H}X'_i\)\(\sum_{i=1}^{e_H}X_i+\sum_{i=1}^{e_G-e_H}X''_i\)\]\\
 &=\p\(G'\cap G''\in  \mathbb{G}(n,p_n)\)\(\E\[\(\sum_{i=1}^{e_H}X_i\)^2\]+\(2e_H(e_G-e_H)+(e_G-e_H)^2\)(\E[X])^2\)\\
 &=p_n^{2e_G-e_H}\(e_H\E\[X^2\]+\(e_G^2-e_H\)\E[X])^2\)\\
 &=p_n^{2e_G-e_H}
  ( e_H\text{Var}[X]+e_G^2(\E[X])^2 ).
\end{align*}
}
Hence we get
\begin{align*}
  \E[S_{G'}S_{G''}]-\E[S_{G'}]\E[S_{G''}]&=p_n^{2e_G-e_H}
  ( e_H\text{Var}[X]+e_G^2(\E[X])^2 )
  -p_n^{2e_G} e_G^2(\E[X])^2
  \\
  &=p_n^{2e_G-e_H}
  (
  e_H\text{Var}[X]+e_G^2(1-p_n^{e_H})(\E[X])^2)
  \\
&\approx p_n^{2e_G-e_H} ( \Var [X]+(1-p_n)(\E[X])^2 ), 
\end{align*}
and consequently
\begin{align*}
\Var\big[W^G_n\big]& \approx
\sum_{\substack{H\subset G\\e_H\geq1}}\sum_{G'\cap G'' \sim H
  \atop
  G', G''\sim G
}
p_n^{2e_G-e_H} ( \Var [X]+(1-p_n)(\E[X])^2 )
\\
&\approx
( \Var [X]+(1-p_n)(\E[X])^2 )\sum_{\substack{H\subset G\\e_H\geq1}}n^{2v_G-v_H}
p_n^{2e_G-e_H} \\
&\approx \big(
    \Var[X]+(1-p_n)(\E[X])^2
    \big)
    \max_{\substack{H\subset G \\ e_H\geq1}} n^{2v_G-v_H}p_n^{2e_G-e_H},
\end{align*}
 as required.
\end{Proof}
Next, we  show in Lemma  \ref{lem:W=sumI} that the combined
weights $W^G_n$ of subgraphs
can be written as a sum of multiple stochastic integrals 
using Proposition~\ref{prop:sum=I}.
This allows us to apply Theorem  \ref{thm:Var<Sum} to
obtain normal approximation in Wasserstein distance for $\widetilde{W}^G_n$,  which is presented in Theorem \ref{thm:main}.
In the sequel we number all possible edges of 
the complete graph $K_n$ from $1$ to $n(n-1)/2$,
and we denote by $E_G \subset \inte^{e_G}$ 
the set of sequences of edges that create a graph isomorphic to $G$,
i.e. a sequence $( e_{k_1} ,\ldots , e_{k_{e_G}} )$ belongs to $E_G$
if and only if the graph created by edges $e_{k_1},\ldots ,e_{k_{e_G}}$
is isomorphic to $G$. Before stating the lemma, let us define the operator $\Psi_{t_i}$
\begin{align}\label{eq:Psi}
  \Psi_{t_i}f(t_1,\ldots , t_n):=f(t_1,\ldots , t_n)-
  \frac{1}{2}
  \int_{2\lfloor t_i/2\rfloor}^{2\lfloor t_i/2\rfloor+2}f(t_1,\ldots ,t_{i-1},s,t_{i+1},\ldots ,t_n)ds,
\end{align}
which arises naturally when reprezenting any multiple stochastic integral $I_n(f_n)$ as $I_n(\bar f_n)$ with $\bar f_n$ satisfying \eqref{ass:int0}, see Proposition \ref{prop:Ifoverline}. 
\begin{lemma}
  \label{lem:W=sumI}
  We have the identity in distribution 
  \begin{equation}
    \label{wg2}
    W^G_n\stackrel{d}{=}\sum_{r=0}^{e_G}  I_k\(\bar{h}_k\),
    \end{equation} 
where
\begin{equation}
  \label{hk} 
\bar{h}_k(t_1,\ldots , t_k ) := \Psi_{t_1} \cdots \Psi_{t_k}g_k \(t_1-2\lfloor t_1/2\rfloor,\ldots ,t_k-2\lfloor t_k/2\rfloor\)\sum_{a\in\N^{e_{\scaleto{G}{3pt}}-k}}\mathbf1_{E_G}\(a,\lfloor t_1/2\rfloor,\ldots ,\lfloor t_k/2\rfloor\), 
\end{equation} 
 and the function $g_k:(0,2)^k\rightarrow \R$ is given by
 \begin{align}
   \label{isgv} 
g_k(t_1,\ldots , t_k)&=
\frac{(p_n/2)^{{e_G}-k}}{( e_G - k)! k!} 
\mathbf1_{(0,2p_n)^k}\(t_1,\ldots , t_k\)\(({e_G}-k)\E[X]+\sum_{i=1}^kF_X^{-1}\(\frac{t_i}{2p_n}\)\),
\end{align}
 where $F_X^{-1}$ is the generalized inverse of the distribution function $F_X$ of $X$.
\end{lemma}
\begin{Proof} 
First, we note that 
\begin{align*}
  W^G_n&\stackrel{d}{=}\frac1{e_G!}\sum_{ k_1\neq\cdots \neq k_{e_G} \geq 0}\mathbf1_{E_G}(k_1, \ldots ,k_{e_G})\mathbf1_{(0,2p_n)^{e_G}}\(U_{k_1}+1,\ldots ,U_{k_{e_G}}+1\)\\
  &\ \ \ \ \ \ \ \ \ \ \ \ \ \ \ \ \ \ \ \ \ \ \ \times\(F^{-1}\(\frac{U_{k_1}+1}{2p_n}\)+
  \cdots
  +F^{-1}\(\frac{U_{k_{e_G}}+1}{2p_n}\)\)\\
  &=\frac1{e_G!}\sum_{ k_1\neq\cdots \neq k_{e_G} \geq 0}
  h_{e_G} (2k_1+1+U_{k_1},\ldots ,2k_{e_G}+1+U_{k_{e_G}} ),
  \end{align*}
 where 
 \begin{eqnarray*}
   h_{e_G}(t_1,\ldots ,t_{e_G})&=&
      \mathbf1_{ E_G}\(\lfloor t_1/2\rfloor,\ldots ,\lfloor t_{e_G}/2\rfloor\)\mathbf1_{(0,2p_n)^{e_G}}\(t_1-2\lfloor t_1/2\rfloor,\ldots ,t_{e_G}-2\lfloor t_{e_G}/2\rfloor\)\\
& & \times\(F^{-1}_X\(\frac{t_1-2\lfloor t_1/2\rfloor}{2p_n}\)+\cdots +F^{-1}_X\(\frac{t_{e_G}-2\lfloor t_{e_G}/2\rfloor}{2p_n}\)\),
\end{eqnarray*}
 and by Proposition \ref{prop:sum=I},
 the relation \eqref{wg2} holds with 
\begin{align*}
h_k(t_1,\ldots , t_k)&:=g_k\(t_1-2\lfloor t_{1}/2\rfloor,\ldots ,t_k-2\lfloor t_k/2\rfloor\)\sum_{a\in\N^{e_{\scaleto{G}{3pt}}-k}}\mathbf1_{E_G}\(a,\lfloor t_1/2\rfloor,\ldots ,\lfloor t_k/2\rfloor\), 
\end{align*}
 where $g_k:(0,2)^k\rightarrow \R$ is given by \eqref{isgv}. 
Finally, in case the functions $h_k$ may
not satisfy the condition \eqref{ass:int0},
we can use Proposition \ref{prop:Ifoverline} to obtain
\eqref{wg2} with
\begin{align*}
&\bar{h}_k(t_1,\ldots , t_k )\\
  &= \Psi_{t_1} \cdots \Psi_{t_k}
  \( g_k \(t_1-2\lfloor t_1/2\rfloor,\ldots ,t_k-2\lfloor t_k/2\rfloor\)\sum_{a\in\N^{e_{\scaleto{G}{3pt}}-k}}\mathbf1_{E_G}\(a,\lfloor t_1/2\rfloor,\ldots ,\lfloor t_k/2\rfloor\)\)
  \\
&= \Psi_{t_1} \cdots \Psi_{t_k}g_k \(t_1-2\lfloor t_1/2\rfloor,\ldots ,t_k-2\lfloor t_k/2\rfloor\)\sum_{a\in\N^{e_{\scaleto{G}{3pt}}-k}}\mathbf1_{E_G}\(a,\lfloor t_1/2\rfloor,\ldots ,\lfloor t_k/2\rfloor\),
\end{align*}
where the last equality follows from the fact that the sum appearing above is constant for $(t_1,\ldots ,t_k)\in (2m_1,2m_1+2)\times\ldots\times(2m_k,2m_k+2)$, $m_1,\ldots,m_k\in\N$. The proof is complete.
\end{Proof}
We can now pass to the main result in this section.
\begin{theorem}\label{thm:main}
 Let  $G$ be a graph without isolated vertices.
 The renormalized weight $\widetilde{W}^G_n$
 of graphs in $\mathbb{G}(n,p_n)$ that are isomorphic to
 $G$ satisfies 
\begin{equation} 
  \label{eq}
  d_W \big(\widetilde{W}^G_n,\mathcal{N}\big) \lesssim
  \frac{\sqrt{\E[(X-\E[X])^4]}+(1-p_n)(\E[X])^2}{\Var [X]+(1-p_n)(\E[X])^2}
  \((1-p_n)\min_{\substack{ H\subset G\\e_H\geq1}} n^{v_H}p_n^{e_H} \)^{-1/2}.
\end{equation}
\end{theorem}
\begin{Proof} 
  Without loss of generality we take $p_n=p$ in the proof.
  By Corollary~\ref{thm:Var<Sum}
 we have 
  \begin{align} 
    \nonumber
    & d_W\big(\widetilde{W}^G_n , \mathcal{N}\big) \lesssim
     \\
\nonumber
      &
\frac1{\Var\[W^G_n\]}
\sqrt{\sum_{0\leq l< k\leq {e_G}}\left\| \bar{h}_{k}\star_{k}^l \bar{h}_k \right\|^2_{
    L^2(\R^{k-l}_+)
    }
   +
    \hskip-0.3cm
    \sum_{1\leq l< k\leq {e_G}}\left\| \bar{h}_l \star_{l}^l \bar{h}_k \right\|^2_{
      L^2(\R^{k-l}_+)
      }
    +
  \hskip-0.3cm
  \sum_{1\leq l< k\leq {e_G}}
  \left\| \bar{h}_k \star_{l}^l \bar{h}_k \right\|^2_{
    L^2(\R^{2(k-l)}_+)
    }
   }
    \\
    \label{eq:dFN<sum}
&=: \frac{\sqrt{S_1+S_2+S_3}}{\Var\[W^G_n\]}, 
\end{align} 
  where $\bar h_k$ has been defined in \eqref{hk}. 
  We note that by the equivalence \eqref{eq:VarW}
  of Proposition~\ref{prop:VarW} it suffices to show that
\begin{equation}
  \label{est:SSS}
  S_1+S_2+S_3\lesssim
 \frac{ 
  \E\[\(X-\E\[X\]\)^4\]+(1-p)^2(\E\[X\])^4
 }{1-p}
 \max_{\substack{H\subset G \\ e_H\geq1}} n^{4v_G-3v_H}p^{4e_G-3e_H}, 
\end{equation} 
which follows from \eqref{est:S1} and \eqref{est:S3} 
below.
 Indeed, applying \eqref{eq:VarW} and \eqref{est:SSS} to \eqref{eq:dFN<sum}
 shows that 
\begin{align*}
  d_W\big(\widetilde{W}^G_n , \mathcal{N}\big)\lesssim\ &\frac{
    \sqrt{ 
      \E\[\(X^2-\E\[X\]\)^4\]+(1-p)^2(\E\[X\])^4 
    \max_{\substack{H\subset G\\e_H\geq1}} n^{4v_G-3v_H}p^{4e_G-3e_H}} 
  }{\sqrt{1-p}
    \big(
    \E [X^2]-p(\E[X])^2
    \big)
    \max_{\substack{H\subset G\\e_H\geq1}} n^{2v_G-v_H}p^{2e_G-e_H}}, 
\end{align*}
and after factoring out $n^{4v_G}p^{4e_G}$ in front of the maxima, we
conclude to 
\begin{align*}
  d_W\big(\widetilde{W}^G_n , \mathcal{N}\big)&\lesssim
  \frac{
    \left( 
    \sqrt{\E\[\(X^2-\E\[X\]\)^4\]}+(1-p)(\E\[X\])^2
    \right)
\(\max_{\substack{H\subset G\\e_H\geq1}} n^{-v_H}p^{-e_H} \)^{3/2}
  }{
    \sqrt{1-p}
    \big(
    \E [X^2]-p(\E[X])^2
    \big)
    \max_{\substack{H\subset G\\e_H\geq1}} n^{-v_H}p^{-e_H} }\\
  &=\frac{\sqrt{\E[(X-\E[X])^4]}+(1-p)(\E[X])^2}{
    \E [X^2]-p(\E[X])^2}\((1-p)\min_{\substack{ H\subset G\\e_H\geq1}} n^{v_H}p^{e_H} \)^{-1/2}. 
\end{align*}
$i)$ Estimation of $S_1$. 
 For $0\leq l< k\leq n$ we have 
\begin{align*}
  &\left\| \bar{h}_k \star_{k}^l \bar{h}_k \right\|^2_{
    L^2(\R_+^{k-l})
    }
  =
  \frac{1}{2^{2l}}
  \int_{\R^{k-l}_+}\(\int_{\R^l_+}\( \bar{h}_k (x_1,\ldots , x_k)\)^2dx_1\cdots dx_l \)^2dx_{l+1}\cdots dx_k
  \\
  &=
  \frac{1}{2^{2l}}
  \int_{\R^{k-l}_+}\left( \sum_{{b}\in\N^l}\(\sum_{a\in\N^{e_{\scaleto{G}{3pt}}-k}}\mathbf1_{E_G}\(a,{b},\lfloor x_{l+1}/2\rfloor,\ldots ,\lfloor x_k/2\rfloor\)\)^2
  \right.
  \\
  &
  \left.
  \quad \int_{(0,2)^l}\Big(
    \Psi_{x_1}\cdots \Psi_{x_k}g_k
    \(x_1,\ldots , x_l,x_{l+1}-2\lfloor x_{l+1}/2\rfloor,\ldots ,x_k-2\lfloor x_k/2\rfloor\)\Big)^2dx_1\cdots dx_l 
    \right)^2dx_{l+1}\cdots dx_k\\
    &=
    \frac{1}{2^{2l}}
    \sum_{{c}\in\N^{k-l}}
  \left( \sum_{{b}\in\N^l}\(\sum_{a\in\N^{e_{\scaleto{G}{3pt}}-k}}\mathbf1_{E_G}\(a,{b},{c}\)\)^2
  \right)^2
  \\
    & \quad 
    \times\int_{(0,2)^{k-l}}\left(
  \int_{(0,2)^l}\Big(
  \Psi_{x_1}\cdots \Psi_{x_k}g_k \(x_1,\ldots , x_k\)\Big)^2dx_1\cdots dx_l
  \right)^2dx_{l+1}\cdots dx_k.
\end{align*}
 Combining the equivalence 
$$ 
  \sum_{{c}\in\N^{k-l}}\left(
  \sum_{{b}\in\N^l}\(\sum_{a\in\N^{e_{\scaleto{G}{3pt}}-k}}\mathbf1_{E_G}\(a,{b},{c}\)\)^2\right)^2
 \approx 
             \max_{\substack{K\subset H\subset G
                 \\ e_K=k-l, \ \! e_H=k}}
             n^{4v_G-2v_H-v_K},
             $$ 
             see the proof of Theorem~4.2 in \cite{PS2},
             with \eqref{eq:Psi2} in Lemma~\ref{l1}, we get 
\begin{align*}
  \left\| \bar{h}_k \star_{k}^l \bar{h}_k \right\|^2_{
    L^2(\R_+^{k-l})
  }\lesssim& \big(
  \E\[\(X-\E\[X\]\)^4\]+(1-p)^2(\E\[X\])^4\big)
  \\
  &\times\max_{\substack{K\subset H\subset G \\
      e_K=k-l, \ \! e_H=k}}
 n^{4v_G-2v_H-v_K}p^{4e_G-2e_H-e_K}(1-p)^{2e_H-e_K-2},
\end{align*}
and consequently 
\begin{eqnarray} 
  \nonumber 
  \lefteqn{
    S_1=\sum_{0\leq l< k\leq n}\left\| \bar{h}_k \star_{k}^l \bar{h}_k \right\|^2_{
    L^2(\R_+^{k-l})
    }
  }
  \\
  \nonumber
  &\lesssim & \big(
  \E\[\(X-\E\[X\]\)^4\]+(1-p)^2(\E\[X\])^4\big)
\max_{\substack{K\subset H\subset G \\ e_K\geq1}}
  n^{4v_G-2v_H-v_K}p^{4e_G-2e_H-e_K}(1-p)^{2e_H-e_K-2}
  \\
  \label{est:S1}
  & \lesssim & 
  \frac{\E\[\(X-\E\[X\]\)^4\]+(1-p)^2(\E\[X\])^4}{1-p} 
  \max_{\substack{H\subset G \\ e_H\geq1}}
n^{4v_G-3v_H}p^{4e_G-3e_H}, 
\end{eqnarray} 
as in the proof of Theorem~4.2 in \cite{PS2}. 
\\ 
             \noindent
$ii)$ Estimation of $S_2$. Similarly, for $1\leq l< k\leq n$ we have 
\begin{align*}
  &\left\| \bar{h}_l \star_{l}^l \bar{h}_k \right\|^2_{
    L^2(\R_+^{k-l})
    }
  =
  \frac{1}{2^{2l}}
  \int_{\R^{k-l}_+}\(\int_{\R^l_+} \bar{h}_l (x_1,\ldots , x_l) \bar{h}_k(x_1,\ldots  ,x_k)dx_1\cdots dx_l\)^2dx_{l+1}\cdots dx_k\\
  &=
  \frac{1}{2^{2l}}
  \int_{\R^{k-l}_+}\Bigg(
  \sum_{b\in\N^l}\(\sum_{a\in\N^{e_{\scaleto{G}{3pt}}-l}}\mathbf1_{E_G}\(a,b\)\sum_{a'\in\N^{e_{\scaleto{G}{3pt}}-k}}\mathbf1_{E_G}\(a',b,\lfloor x_{l+1}/2\rfloor,\ldots ,\lfloor x_k/2\rfloor\)\)\\
    &\quad \times \int_{(0,2)^l}
    \Psi_{x_1}\cdots \Psi_{x_l}g_l (x_1, \ldots , x_l)
    \Psi_{x_1}\cdots \Psi_{x_k}g_k
    \(x,x_{l+1}-2\lfloor x_{l+1}/2\rfloor,\ldots ,x_k-2\lfloor x_k/2\rfloor\)dx_1\cdots dx_l \Bigg)^2
    \\
    &\quad dx_{l+1}\cdots dx_k
  \\
  &=
  \frac{1}{2^{2l}}
  \sum_{c\in\N^{k-l}}\Bigg(
  \sum_{b'\in\N^l}\(\sum_{a\in\N^{e_{\scaleto{G}{3pt}}-l}}\mathbf1_{E_G}\(a,b\)\sum_{a'\in\N^{e_{\scaleto{G}{3pt}}-k}}\mathbf1_{E_G}\(a',b,c\)\)\Bigg)^2\\
  &\quad \times\int_{(0,2)^{k-l}}\Bigg(
  \int_{(0,2)^l}
    \Psi_{x_1}\cdots \Psi_{x_l}g_l(x_1,\ldots , x_l)
    \Psi_{x_1}\cdots \Psi_{x_k}g_k(x_1,\ldots , x_k)dx_1\cdots dx_l\Bigg)^2dx_{l+1}\cdots dx_k.
\end{align*}
By the Cauchy-Schwarz inequality and the formula \eqref{eq:Psi1} in Lemma~\ref{l1},
we get 
\begin{align*}
  &\int_{(0,2)^{k-l}}\Bigg( \int_{(0,2)^l}
    \Psi_{x_1}\cdots \Psi_{x_l}g_l (x_1,\ldots , x_l)
    \Psi_{x_1}\cdots \Psi_{x_k}g_k\(x_1,\ldots , x_k\)dx_1\cdots dx_l\Bigg)^2
    dx_{l+1}\cdots dx_k\\
  &\leq \int_{(0,2)^l}
  \Psi_{x_1}\cdots \Psi_{x_l}g_l^2(x_1\cdots x_l)dx_1\cdots dx_l
  \int_{(0,2)^k}
  \Psi_{x_1}\cdots \Psi_{x_k} g_k^2\(x_{l+1}\cdots x_k\)dx_{l+1}\cdots dx_k
  \\
&\lesssim p^{4{e_G}-k-l}(1-p)^{k+l-2}\big( \E\big[(X^2-\E\[X\])^2\big]+(1-p)(\E\[X\])^2\big)^2\\
&\lesssim \frac{p^{4{e_G}-k-l}}{1-p}\big( \E\big[(X-\E\[X\])^4\big]+(1-p)^2(\E\[X\])^4\big).
\end{align*}
Furthermore, we have 
\begin{equation} 
  \label{dssds2}
    \sum_{c\in\N^{k-l}}\Bigg(
    \sum_{a'\in\N^l}\(\sum_{a'\in\N^{e_{\scaleto{G}{3pt}}-l}}\mathbf1_{E_G}\(a,b\)\sum_{a'\in\N^{e_{\scaleto{G}{3pt}}-k}}\mathbf1_{E_G}\(a',b,c\)\)\Bigg)^2
\lesssim \max_{\substack{K\subset H' \subset G\\e_K=k-l, \ \! e_{H'} =l }} n^{4v_G-2v_{H'} -v_K}, 
\end{equation} 
see the proof of Theorem~4.2 in \cite{PS2}, thus
\begin{align*}
  \left\| \bar{h}_l \star_{l}^l \bar{h}_k \right\|^2_{
   L^2(\R_+^{k-l})
  }\lesssim&
  \frac{\E\big[(X-\E\[X\])^4\big]+(1-p)^2(\E\[X\])^4}{1-p}
\max_{\substack{K\subset H'\subset G \\ e_K=k-l, \ \! e_{H'}=k}}
 n^{4v_G-2v_{H'}-v_K}p^{4e_G-2e_{H'}-e_K},
\end{align*}
 from which it follows by that 
\begin{equation} 
 \label{est:S3}
 S_2 = \sum_{1\leq l< k\leq n}\left\| \bar{h}_l \star_{l}^l \bar{h}_k \right\|^2_{
    L^2(\R_+^{k-l}) }
 \lesssim 
    \frac{E\big[(X-\E\[X\])^4\big]+(1-p)^2(\E\[X\])^4}{1-p}
 \max_{\substack{H\subset G \\ e_H\geq1}}
 n^{4v_G-3v_H}p^{4e_G-3e_H}, 
\end{equation} 
as in the proof of Theorem~4.2 in \cite{PS2}. 
\\ 
  $iii)$ Estimation of $S_3$. 
For $1\leq l< k\leq n$ we have
\begin{align*}
  &\left\| \bar{h}_k \star_{l}^l \bar{h}_k \right\|^2_{
    L^2(\R_+^{2(k-l)})
    }\\
  &=
  \frac{1}{2^{2l}}
  \int_{\R^{k-l}_+}\int_{\R^{k-l}_+}\(\int_{\R^l_+}
  \bar{h}_k (x_1,\ldots , x_k)
  \bar{h}_k (x_1,\ldots , x_l,z_1,\ldots , z_{k-l}
  )dx_1\cdots dx_l\)^2dx_{l+1}\cdots dx_k dz_1\cdots dz_{k-l}\\
  &=
  \frac{1}{2^{2l}}
  \int_{\R^{k-l}_+}\int_{\R^{k-l}_+}\Bigg(
  \sum_{b\in\N^l}
  \sum_{a,a'\in\N^{e_{\scaleto{G}{3pt}}-k}}
  \mathbf1_{E_G}\(a,b,\lfloor x_{l+1}/2\rfloor,\ldots ,\lfloor x_k/2\rfloor\)
  \mathbf1_{E_G}\(a',b,\lfloor z_1/2\rfloor,\ldots ,\lfloor z_{k-l}/2\rfloor\)
  \\
    &\quad \int_{(0,2)^l}
  \Psi_{x_1}\cdots \Psi_{x_k}g_k\(x,x_{l+1}-2\lfloor x_{l+1}/2\rfloor,\ldots ,x_k-2\lfloor x_k/2\rfloor\)
  \\
&\quad \Psi_{x_1}\cdots \Psi_{x_l}\Psi_{z_1}\cdots \Psi_{z_{k-l}}g_k\(x,x_{l+1}-2\lfloor x_{l+1}/2\rfloor,\ldots ,z_{k-l}-2\lfloor z_{k-l}/2\rfloor\)dx\Bigg)^2dx_{l+1}\cdots dx_k dz_1\cdots dz_{k-l}.
\end{align*}
Then, applying the Cauchy-Schwarz inequality to the inner integral,  we get
\begin{eqnarray*}
    \left\| \bar{h}_k \star_{l}^l \bar{h}_k \right\|^2_{
    L^2(\R_+^{2(k-l)})
    }
 & = & 
  \frac{1}{2^{2l}}
  \sum_{c,c'\in\N^{k-l}}\Bigg(
  \sum_{b\in\N^l}\(\sum_{a\in\N^{e_{\scaleto{G}{3pt}}-k}}\mathbf1_{E_G}\(a,b,c\)\)\(\sum_{a'\in\N^{e_{\scaleto{G}{3pt}}-k}}\mathbf1_{E_G}\(a',b,c'\)\)\Bigg)^2
  \\
& & \times\(\int_{(0,2)^k}\Big(\Psi_{x_1}\cdots \Psi_{x_k}g_k\(x_1,\ldots , x_k\)\Big)^2dx_1\cdots dx_k\)^2.
\end{eqnarray*}
Since $k\geq1$, the formula \eqref{eq:Psi1} in Lemma~\ref{l1} gives us
\begin{eqnarray*}
  \lefteqn{
    \(\int_{(0,2)^k}\Big(\Psi_{x_1}\cdots \Psi_{x_k}g_k\(x_1,\ldots , x_k\)\Big)^2dx_1\cdots dx_k
  \)^2
  }
  \\
   & \lesssim & p^{4{e_G}-2k}(1-p)^{2k-2}
    \big(\E\big[
      \(X^2-\E\[X\]\)^2\big]+(1-p)(\E\[X\])^2\big)^2
  \\
  &\lesssim & 
  \displaystyle
  \frac{p^{4{e_G}-2k}}{1-p}
  \big( \E\big[\(X-\E\[X\]\)^4\big]+(1-p)^2(\E\[X\])^4\big).
\end{eqnarray*}
Furthermore, we have
\begin{align*}
\sum_{c,c'\in\N^{k-l}}&\Bigg( \sum_{b\in\N^l}\(\sum_{a\in\N^{e_{\scaleto{G}{3pt}}-k}}\mathbf1_{E_G}\(a,b,c\)\)\(\sum_{a'\in\N^{e_{\scaleto{G}{3pt}}-k}}\mathbf1_{E_G}\(a',b,c'\)\)\Bigg)^2\\[8pt]
 & \lesssim \max_{\substack{K,H,L\subset G\\e_K =k-l, \ \! e_H =l,\ e_L=k}}
  n^{4v_G-v_K-v_H-v_L},
\end{align*} 
see the proof of Theorem~4.2 in \cite{PS2}, 
from which it follows
$$ 
 S_3 =
  \sum_{1\leq l< k\leq e_G}
\left\| \bar{h}_k \star_{l}^l \bar{h}_k \right\|^2_{
  L^2(\R_+^{2(k-l)})
}
\lesssim
\frac{\E\big[(X-\E\[X\])^4\big]+(1-p)^2(\E\[X\])^4}{1-p}
\max_{\substack{H\subset G\\e_H \geq1}}\ \! n^{4v_G-3v_H }p^{4e_G-3e_H },
$$
as in the proof of Theorem~4.2 in \cite{PS2},
which concludes the proof by \eqref{eq:VarW} and \eqref{eq:dFN<sum}.
\end{Proof}
We note that the bound \eqref{eq} implies 
$$
  d_W \big(\widetilde{W}^G_n,\mathcal{N}\big) \lesssim
  \left(
  \frac{(\E[X])^2}{\Var [X]}
  +
  \sqrt{\kappa_X}
  \right)
  \((1-p_n)\min_{\substack{ H\subset G\\e_H\geq1}} n^{v_H}p_n^{e_H} \)^{-1/2},
  $$
  where
  $\E[X]/\sqrt{\Var [X]}$ is the standardized first moment of $X$ and 
  $$
  \kappa_X : = \frac{\E[(X-\E[X])^4]}{ ( \Var [X] )^2 }
  $$
  is the kurtosis of $X$.
  \\
   
In the next corollary
 we note that Theorem~\ref{thm:main} simplifies 
 if we narrow our attention to $p_n$ depending of the complete
 graph size $n$ and close to $0$ or to $1$.
\begin{corollary}
   \label{c0}
Let  $G$ be a graph without separated vertices. For $p_n<c<1$, $n\geq 1$, we have
\begin{align*}
d_W \big(\widetilde{W}^G_n,\mathcal{N}\big)&\lesssim\frac{\sqrt{\E[X^4]}}{\E[X^2]}\((1-p_n)
  \min_{\substack{ H\subset G\\e_H\geq1}} n^{v_H}p_n^{e_H} \)^{-1/2}.
\end{align*}
  On the other hand, for  $p_n>c>0$, $n\geq 1$, it holds
 \begin{equation}
   \label{dk}
d_W \big(\widetilde{W}^G_n,\mathcal{N}\big)\lesssim
\frac{\sqrt{\E[X^4]}}{
  n\sqrt{1-p_n} \Var [X] }.
 \end{equation} 
\end{corollary}
 Furthermore, it turns out that the mininum appearing in Theorem \ref{thm:main} and above for a wide class of graphs satisfying a certain balance condition.
Precisely, let us consider the class $\mathcal B$ of all
graphs with at least three vertices, and such that 
\begin{equation}\label{eq:balance}
\max_{\substack{H\subset G\\ {v_H\geq3}}}\frac{e_H-1}{v_H-2}=\frac{e_G-1}{v_G-2}, 
\end{equation} 
 as introduced in \cite{PS2}.
 It has been shown there that a graph with at least 3 vertices and at least
 one edge belongs to $\mathcal B$ if and only if
  for any $p\in(0,1)$ {and $n\geq v_G$} we have 
  $$\min_{\substack{ H\subset G\\e_H\geq1}}n^{v_H}p^{e_H}=\min\{n^2p, n^{v_G}p^{e_G}\}
.$$
  An application of this fact to Corollary \ref{c0} yields
  the following result. 
\begin{prop}\label{prop:dKforGinB}
For $G\in \mathcal B$ without isolated vertices and $c\in(0,1)$ we have
\begin{align*}
  d_W \big(\widetilde{W}^G_n,\mathcal{N}\big)\lesssim
    \left\{\begin{array}{ll}
  \displaystyle
 \frac{\sqrt{\E[X^4]}}{n\sqrt{1-p_n}\Var [X]} &\mbox{if } \ \displaystyle
  0<c< p_n, 
 \\
 \\
  \displaystyle
  \frac{\sqrt{\E[X^4]}}{n\sqrt{p_n}\E [X^2]}
 &
  \displaystyle
  \mbox{if } \ n^{-(v_g-2)/(e_G-1)} < p_n \leq c, 
  \\
  \\
  \displaystyle
  \frac{\sqrt{\E[X^4]}}{n^{v_G/2}p_n^{e_G/2}\E [X^2]}
  &
  \displaystyle
  \mbox{if } \ 0 < p_n \leq n^{-(v_G-2)/(e_G-1)}.
  \end{array}\right.
\end{align*}
\end{prop}
\noindent
 The following Corollaries~\ref{cor1.3}-\ref{cor1.5} of Proposition \ref{prop:dKforGinB}
can be proved similarly to their counterparts Corollaries~4.8-4.10 in \cite{PS2}. 
The next Corollary~\ref{cor1.3} deals with cycle graphs with $r$ vertices, $r\geq3$,
and in particular with triangles when $r=3$.
\begin{corollary}
\label{cor1.3}
Let $G$ be a cycle graph with $r$ vertices, $r\geq3$, and
$c\in (0,1)$. We have 
\begin{align*}
  d_W \big(\widetilde{W}^G_n,\mathcal{N}\big)\lesssim
    \left\{\begin{array}{ll}
  \displaystyle
 \frac{\sqrt{\E[X^4]}}{n\sqrt{1-p_n}\Var [X]} &\mbox{if } \ \displaystyle
  0<c< p_n, 
 \\
 \\
  \displaystyle
  \frac{\sqrt{\E[X^4]}}{n\sqrt{p_n}\E [X^2]}
 &
  \displaystyle
  \mbox{if } \ n^{-(r-2)/(r-1)} < p_n \leq c, 
  \\
  \\
  \displaystyle
  \frac{\sqrt{\E[X^4]}}{(np_n)^{r/2}\E [X^2]}
  &
  \displaystyle
  \mbox{if } \ 0 < p_n \leq n^{-(r-2)/(r-1)}.
  \end{array}\right.
\end{align*}
\end{corollary}
In the case of complete graphs,
the next corollary also covers the case of
triangles. 
\begin{corollary}
\label{cor1.4}
 Let $G$ be a complete graph with $r$ vertices,
 $r\geq3$, and $c\in (0,1)$. We have 
 \begin{align*}
  d_W \big(\widetilde{W}^G_n,\mathcal{N} \big)\lesssim
    \left\{\begin{array}{ll}
  \displaystyle
    \frac{\sqrt{\E[X^4]}}{n\sqrt{1-p_n}\Var [X]} 
 &\mbox{if } \ \displaystyle
  c< p_n <1, 
 \\
 \\
    \displaystyle
    \frac{\sqrt{\E[X^4]}}{n\sqrt{p_n}\E[X^2]} 
 &
      \displaystyle
      \mbox{if } \
 n^{-2/(r+1)} < p_n \leq c, 
      \\
      \\
  \displaystyle
  \frac{\sqrt{\E[X^4]}}{ n^{r/2}p_n^{r(r-1)/4} \E[X^2]} 
  &  \displaystyle
 \mbox{if } \ 0 < p_n \leq n^{-2/(r+1)}. 
  \end{array}\right.
\end{align*}
\end{corollary}
 Finally, the last corollary deals with
the important class of graphs which have a tree structure.
\begin{corollary}
  \label{cor1.5}
  Let $G$ be any tree (a connected graph without cycles) with $r$ edges,
  and $c\in (0,1)$.
  We have 
\begin{align*}
  d_W \big(\widetilde{W}^G_n,\mathcal{N} \big)\lesssim
     \left\{\begin{array}{ll}
    \displaystyle
    \frac{\sqrt{\E[X^4]}}{ n\sqrt{1-p_n} \Var [X]}
     &\mbox{if } \  \displaystyle
 c< p_n<1, 
 \\
 \\
    \displaystyle
    \frac{\sqrt{\E[X^4]}}{ n\sqrt{p_n} \E[X^2]} 
     &  \displaystyle
    \mbox{if } \ \frac{1}{n} < p_n \leq c,\\ \\
  \displaystyle
  \frac{\sqrt{\E[X^4]}}{ n^{(r+1)/2}p_n^{r/2} \E[X^2]} 
  &  \displaystyle
   \mbox{if } \ 0 < p_n \leq \frac{1}{n}.
  \end{array}\right.
\end{align*}
\end{corollary}

\section{Appendix}
\label{sec-app} 
In this section we gather a number of technical results,
starting with the following multiplication formula
for multiple stochastic integrals, which involves the $\star$-notation introduced in \eqref{starnotation}. 
 For $f_n \in \widehat{L}^2(\R_+^n)$
and $g_m \in \widehat{L}^2(\R_+^m)$
satisfying \eqref{ass:int0} the following multiplication formula holds: 
\begin{equation}
 \label{eq:mf}
 I_n(f_n) I_m(g_m)=\sum_{k=0}^{m\wedge n}k!{{m}\choose{k}}{{n}\choose{k}}\sum_{i=0}^k{{k}\choose{i}} I_{m+n-k-i}\big(
 f_n\hskip0.1cm \widetilde{\star}_{k}^{i} g_m\big),
\end{equation}
whenever $f_n\star_{k}^{i}g_m\in L^2(\R_+^{m+n-k-i})$ for every $0\leq i \leq k\leq m\wedge n$, see Proposition~5.1 of \cite{PS}.  The next proposition allows us
to bound the $L^2$ norm of $f_n\star g_m$ 
by some simpler expressions, which is used in the proof of Theorem \ref{thm:Var<Sum}. 
\begin{prop}
\label{prop:fg<f^2+g^2}
 Let $f_n\in {L}^2(\R_+^n)$ and $g_m \in {L}^2(\R_+^m)$
 be symmetric functions. 
 For  $0\leq l< k\leq n \wedge m$ we have 
\begin{align}\label{eq:l<k}
  \left\| f_n\star_{k}^l g_m\right\|^2_{{L}^2(\R_+^{m+n-k-l})}2^{2n-2k-1} \leq
  \left\| f_n\star_n^{l+n-k} f_n\right\|^2_{L^2(\R_+^{k-l})}
  + 2^{2m-2k-1}
  \left\| g_m\star_m^{l+m-k} g_m\right\|^2_{L^2(\R_+^{k-l})},
\end{align}
and for $0\leq  k\leq n \wedge m$ we have
\begin{equation} 
\label{eq:l=k}
  \left\| f_n\star_{k}^k g_m\right\|^2_{{L}^2(\R_+^{m+n-2k})}
  \leq  
 2^{2n-4k-1}
 \left\| f_n\star_{n-k}^{n-k} f_n\right\|^2_{L^2(\R_+^{2k})}
 +
 2^{2m-4k-1}
 \left\| g_m\star_{m-k}^{m-k} g_m\right\|^2_{L^2(\R^{2k}_+)}. 
\end{equation}
\end{prop}
\begin{Proof}
 Let $x\in \R^l_+$, $y\in \R^{k-l}_+$, $u\in \R^{n-k}_+$ and $z\in\R^{m-k}_+$.
 H\"older's inequality applied twice gives us
\begin{eqnarray}
  \nonumber
  \lefteqn{
    \left\|
    f_n\star_{k}^l g_m\right\|^2_{L^2(\R_+^{m+n-k-l})}
  }
  \\
  \nonumber
   & = & 
  \frac{1}{2^{2l}}
  \int_{\R^{m-k}_+}\int_{\R^{n-k}_+}\int_{\R^{k-l}_+}
  \(\int_{\R^l_+}f_n(x,y,u)g_m(x,y,z)dx\)^2dydudz
  \\
  \nonumber
  &\leq &
  \frac{1}{2^{2l}}
  \int_{\R^{k-l}_+}
  \int_{\R^{m-k}_+}\int_{\R^{n-k}_+}
  \int_{\R^l_+}f_n^2(x,y,u)dx
  \int_{\R^l_+}g_m^2(x,y,z)dxdudzdy
  \\
  \nonumber
  &\leq & 
    \frac{1}{2^{2l}}
    \sqrt{
      \int_{\R_+^{k-l}}
      \(\int_{\R^{n-k}_+}\int_{\R^l_+}f_n^2(x,y,u)dxdu\)^2dy
      \int_{\R_+^{k-l}}
            \(\int_{\R^{m-k}_+}\int_{\R^l_+}g_m^2(x,y,z)dxdz\)^2dy
  }
  \\
  \nonumber
  &\leq & 
  \frac{1}{2^{2l+1}}
  \int_{\R_+^{k-l}}
  \(\int_{\R^{n-k}_+}\int_{\R^l_+}f_n^2(x,y,u)dxdu\)^2dy
 + 
  \frac{1}{2^{2l+1}}
  \int_{\R_+^{k-l}}
  \(\int_{\R^{m-k}_+}\int_{\R^l_+}g_m^2(x,y,z)dxdz\)^2dy
  \\
  & = &
   \label{eq:l<k2}
  2^{2n-2k-1} 
  \left\| f_n\star_n^{l+n-k} f_n\right\|^2_{L^2(\R_+^{k-l})}
  + 2^{2m-2k-1}
  \left\| g_m\star_m^{l+m-k} g_m\right\|^2_{L^2(\R_+^{k-l})}, 
\end{eqnarray}
where we used the inequality $\sqrt {ab} \leq (a+b)/2$, $a,b\geq0$, 
which proves the first assertion.
Furthermore, for $x,u\in\R^k_+$, $y\in\R^{n-k}_+$ and $z\in\R^{m-k}_+$ 
we get
\begin{align} 
\nonumber 
& \left\|    
f_n\star_{k}^k g_m\right\|^2_{L^2(\R_+^{m+n-2k})}
\\
\nonumber 
 & =  
   \frac{1}{2^{2k}}
   \int_{\R^{n-k}_+}\int_{\R^{m-k}_+}
   \int_{\R_+^k} f_n(u,y)g_m(u,z)du
   \int_{\R_+^k} f_n(x,y)g_m(x,z)dx
   dydz
\\
\nonumber 
& \leq  
   \frac{1}{2^{2k}}
   \int_{\R_+^k}\int_{\R_+^k}
   \(\int_{\R^{n-k}_+}f_n(u,y)f_n(x,y)dy\)\(\int_{\R^{m-k}_+}g_m(u,z)g_m(x,z)dz\)dudx
   \\
\nonumber 
   &\leq 
   \frac{1}{2^{2k+1}}
   \left(
   \int_{\R_+^k}\int_{\R_+^k}\(\int_{\R^{n-k}_+}f_n(u,y)f_n(x,y)dy\)^2dudx
   +
   \int_{\R_+^k} \int_{\R_+^k} \(\int_{\R^{m-k}_+}g_m(u,z)g_m(x,z)dz\)^2dudx
   \right)
   \\
\nonumber 
   & \leq  
 2^{2n-4k-1}
 \left\| f_n\star_{n-k}^{n-k} f_n\right\|^2_{L^2(\R_+^{2k})}
 +
 2^{2m-4k-1}
 \left\| g_m\star_{m-k}^{m-k} g_m\right\|^2_{L^2(\R^{2k}_+)}.
\end{align}
\end{Proof}
The next proposition presents some relationships between second norms involving operators $\nabla$, $L$ and $(-L)^{1/2}$.
\begin{prop}\label{prop:sm} For $F$ such that $LF\in L^2(\Omega)$ we have
\begin{align}\label{eq:sm}\E \left[ \int_0^\infty \( \nabla_t F \)^2 \frac{dt}2\right]=  \E\[ \( (-L)^{1/2} F \)^2\]\leq \E\[(LF)^2\].\end{align}
\end{prop}
\begin{Proof}
Using the chaos decomposition \eqref{eq:decom},
where the sequence of functions $f_n$ in $\widehat{L}^2 (\R_+^n)$, $n\geq 1$, 
satisfies the Condition~\eqref{ass:int0}, and by the isometry relation~\eqref{eq:EI^2}
we have
\begin{align*}
  \E \left[ \int_0^\infty | \nabla_t  F |^2 \frac{dt}2\right] &=
  \sum_{n=1}^\infty
  \E \left[ \int_0^\infty |  \nabla_t  I_n (f_n) |^2 \frac{dt}2\right]
  \\
  &=
  \sum_{n=1}^\infty
  \E \left[ \int_0^\infty |  \nabla_t I_n (f_n) |^2 \frac{dt}2\right]
  \\
  &=
  \sum_{n=1}^\infty
  n^2\E \left[ \int_0^\infty |  I_{n-1}(f_n(t,\cdot)) |^2 \frac{dt}2\right]
  \\
  &=
  \sum_{n=1}^\infty
  n^2(n-1)!\int_0^\infty  \|f_n(t,\cdot)\|^2_{\widehat{L}^2(\R_+^{n-1},dx/2)}\frac{dt}2
  \\
  &=  \sum_{n=1}^\infty {n} \E\[|I_n(f_n)|^2\]
   \\
   &=   \sum_{n=1}^\infty \E\[ | (-L)^{1/2} I_n(f_n)|^2\]
   \\
   &=   \E\big[ \big( (-L)^{-1/2} F \big)^2\big]
   ,
\end{align*}
which is the first part of the assertion. This also implies
$$\E\big[ \( (-L)^{-1/2} F \)^2\big]=\sum_{n=1}^\infty {n} \E\[|I_n(f_n)|^2\]\leq \sum_{n=1}^\infty {n^2} \E\[|I_n(f_n)|^2\]\leq \E [ \( L F \)^2],
$$
which ends the proof.
\end{Proof}
Next, let us recall the definition \eqref{eq:Psi} of the operator $\Psi_{t_i}$
\begin{align*}
  \Psi_{t_i}f(t_1,\ldots , t_n):=f(t_1,\ldots , t_n)-
  \frac{1}{2}
  \int_{2\lfloor t_i/2\rfloor}^{2\lfloor t_i/2\rfloor+2}f(t_1,\ldots ,t_{i-1},s,t_{i+1},\ldots ,t_n)ds,
\end{align*}
$i=1,\ldots , n$,
$t_1,\ldots ,t_n \in \real_+$. The following result is the analog
of the Stroock formula~\cite{stroock}
in our framework.
 
\begin{prop}\label{prop:Ifoverline}
  For every $f_n\in \widehat{L}^2(\R_+^n)$
  there exists a unique $\bar{f}_n\in \widehat{L}^2(\R_+^n)$
  satisfying \eqref{ass:int0} such that $I_n(f_n)=I_n (\bar{f}_n )$, and it is given by
  \begin{align}\label{eq:Ifoverline}
  \bar{f}_n(t_1,\ldots ,t_n)
= 
\Psi_{t_1} \cdots \Psi_{t_n} f_n(t_1,\ldots ,t_n)=\frac1{n!}\nabla_{t_1} \cdots \nabla_{t_n}I_n(f_n).
\end{align}

\end{prop}
\begin{Proof} 
 Uniqusness of $\bar f_n$ follows from the isometry
relation~\eqref{eq:EI^2}. We can also check that the condition~\eqref{ass:int0}
is satisfied by integrating \eqref{eq:Psi} with respect to $t_i \in \real_+$. 
Furthermore,  the equality \eqref{eq:Ifoverline} is clear for $n=1$.
  Assuming that it holds for some $n-1\geq1$, we get
  \begin{align*}
  I_n(f_n)&=\int_0^\infty I_{n-1}\(f_n(t_1,*)\)d(Y_{t_1} - t_1 / 2 )\\
  &=\int_0^\infty I_{n-1}\(\Psi_{t_2} \cdots \Psi_{t_n}f_n(t_1,*)\)d(Y_{t_1} - t_1 / 2 )\\
  &=\int_0^\infty \Psi_{t_1}I_{n-1}\(\Psi_{t_2} \cdots \Psi_{t_n}f_n(t_1,*)\)dY_{t_1} \\
   &=\int_0^\infty I_{n-1}\(\Psi_{t_1} \cdots \Psi_{t_n}f_n(t_1,*)\)dY_{t_1} \\
   &=\int_0^\infty I_{n-1}\(\Psi_{t_1} \cdots \Psi_{t_n}f_n(t_1,*)\)d(Y_{t_1}-t_1/2) \\
   &= I_{n}\(\Psi_{t_1} \cdots \Psi_{t_n}f_n\).  
\end{align*}
Eventually, the latter equality in \eqref{eq:Ifoverline} follows from \eqref{eq:nablaext}.
\end{Proof}
\noindent 
The following Lemma~\ref{l1} is used to bound the 
kernel functions $\bar{h}_k$ appearing in Lemma~\ref{lem:W=sumI}.

\begin{lemma}
  \label{l1}
  The functions $g_k$ defined in \eqref{isgv} 
  satisfy the inequalities 
\begin{align}\nonumber
  &\int_{(0,2)^k} \(\Psi_{x_1}\cdots \Psi_{x_k}g_k\(x_1,\ldots , x_k\)\)^2dx_1\cdots dx_k
  \\\label{eq:Psi1}
  &\hspace{10mm}\lesssim p^{2{e_G}-k}(1-p)^{k-1}\big(
  \E\big[\(X^2-\E\[X\]\)^2\big]+(1-p)(\E\[X\])^2\big)
\end{align}
and
\begin{align}\nonumber
  &\int_{(0,2)^{k-l}}\left(
  \int_{(0,2)^l}\big( \Psi_{x_1}\cdots \Psi_{x_k}g_k\(x_1,\ldots , x_k\)\big)^2dx_1\cdots dx_l
  \right)^2dx_{l+1}\cdots dx_k\\
  \label{eq:Psi2}
  &\hspace{10mm}\lesssim p^{4{e_G}-3k+l}(1-p)^{k+l-2}
  \big(
  \E\[\(X-\E\[X\]\)^4\]+(1-p)^2(\E\[X\])^4
  \big),
\end{align}
$0 \leq l \leq k \leq n$,
where the operator $\Psi_x$ is defined in \eqref{eq:Psi}. 
\end{lemma}
\begin{Proof}
  We decompose $g_k(x_1,\ldots , x_k) $ as
  $$
  g_k(x_1,\ldots , x_k) =\sum_{i=0}^kg_k^{(i)}(x_1,\ldots , x_k) 
  $$
  where 
 $$
 g_k^{(0)}(x_1,\ldots , x_k) :=
 \frac{(p/2)^{{e_G}-k}}{( e_G-k)!k!}
 ({e_G}-k)\E[X]
 \mathbf1_{(0,2p)^k}\(x_1,\ldots , x_k\), 
 $$
 and
$$ 
g_k^{(i)}(x_1,\ldots , x_k) :=
\frac{(p/2)^{{e_G}-k}}{( e_G-k)!k!}
\mathbf1_{(0,2p)^k}\(x_1,\ldots , x_k\)
F^{-1}_X\(\frac{x_i}{2p}\),\qquad 1\leq i\leq k.
$$
Next, for $1\leq i\leq k$ we have 
\begin{align*}
  &\Psi_{x_1}\cdots \Psi_{x_k}g_k^{(i)}(x_1,\ldots , x_k)
  =
  \frac{(p/2)^{{e_G}-k}}{(e_G-k)!k!}
  \left(
  \mathbf1_{(0,2p)}\(x_i\)F^{-1}_X\(\frac{x_i}{2p}\)-p\E\[X\]
  \right)
  \prod_{\substack{1\leq j\leq k\\j\neq i}}\(\mathbf1_{(0,2p)}\(x_j\)-p\).
\end{align*}
 Thus we have 
\begin{eqnarray*}
  \lefteqn{
    \int_{(0,2)^k}\big(\Psi_{x_1}\cdots \Psi_{x_k}g^{(i)}_k\(x_1,\ldots , x_k \)\big)^2
    dx_1\cdots dx_k
  }
  \\
  &=&
\frac{2^kp^{2{e_G}-k-1}(1-p)^{k-1}
}{((e_G-k)!k!)^22^{2{e_G}-2k}}
  \big(
  p\E\big[ (X-p\E\[X\] )^2\big]+(1-p)(p\E\[X\])^2
  \big)
  \\
  &\lesssim & p^{2{e_G}-k}(1-p)^{k-1}\big(
  \E [X^2]-p(\E\[X\])^2\big),
  \qquad
  1\leq i\leq k,
\end{eqnarray*}
 and similarly for $g_k^{(0)}(x_1,\ldots , x_k)$,
which gives 
\begin{eqnarray*}
  \lefteqn{
    \! \! \! \! \! \! \! \! \! \! \! \! \! \! \! \! \! \! \! \! \! \! \! \! \! \! \! \! \!
       \int_{(0,2)^k}\(\Psi_{x_1}\cdots \Psi_{x_k}g_k\(x_1,\ldots , x_k\)\)^2
  dx_1\cdots dx_k
  }
  \\
  & 
  \lesssim & \sum_{i=0}^k\int_{(0,2)^k}\big(
  \Psi_{x_1}\cdots \Psi_{x_k}g^{(i)}_k(x_1,\ldots , x_k )\big)^2dx_1\cdots dx_k
\\
  &\lesssim & p^{2{e_G}-k}(1-p)^{k-1}\big(
  \E [ X^2 ] -p (\E\[X\])^2 \big),
\end{eqnarray*}
as required. In order to prove \eqref{eq:Psi2}, 
we proceed similarly and get
\begin{eqnarray*}
  \lefteqn{ 
    \! \! \! \! \! \! \! \! \! \! \! \! \! \! \! \! \! \! \! \! \!
    \int_{(0,2)^{k-l}}\left(
  \int_{(0,2)^l}\Big(\Psi_{x_1}\cdots \Psi_{x_k}g^{(i)}_k\(x_1,\ldots ,x_k\)\Big)^2dx_1\cdots dx_l
  \right)^2dx_{l+1}\cdots dx_k
  }
  \\
  &\lesssim & p^{4{e_G}-3k+l}(1-p)^{k+l-2}
  \big( \E [X^2]-p(\E\[X\])^2\big)^2
\end{eqnarray*}
for $1\leq i\leq l$, and 
\begin{eqnarray*}
  \lefteqn{
    \! \! \! \! \! \! \! \! \! \! \! \! \! \! \! \! \! \! \! \! \!
    \int_{(0,2)^{k-l}}
  \left(
  \int_{(0,2)^l}\Big(\Psi_{x_1}\cdots \Psi_{x_k}g^{(i)}_k\(x_1,\ldots , x_k\)\Big)^2dx_1\cdots dx_l
  \right)^2dx_{l+1}\cdots dx_k
  }
  \\[9pt]
  &\lesssim & p^{4{e_G}-3k+l}(1-p)^{k+l-1}\big(
    \E\[\(X-\E\[X\]\)^4\]+(1-p)(\E\[X\])^4\big)
\end{eqnarray*}
for $l< i\leq k$. Hence, by the Cauchy-Schwarz inequality we get 
\begin{eqnarray*}
  \lefteqn{ 
  \int_{(0,2)^{k-l}}\left(
  \int_{(0,2)^l}\big(
  \Psi_{x_1}\cdots \Psi_{x_k}g_k
  (x_1,\ldots , x_k )\big)^2dx_1\cdots dx_l\right)^2dx_{l+1}\cdots dx_k
  }
  \\
  &\lesssim & \sum_{i=0}^k\int_{(0,2)^{k-l}}\left(
  \int_{(0,2)^l}\big(\Psi_{x_1}\cdots \Psi_{x_k}g^{(i)}_k
  (x_1,\ldots , x_k )\big)^2dx_1\cdots dx_l
  \right)^2dx_{l+1}\cdots dx_k\\[9pt]
  &\lesssim & p^{4{e_G}-3k+l}(1-p)^{k+l-2}\big(
  \E\[\(X-\E\[X\]\)^4\]+(1-p)^2(\E\[X\])^4\big),
\end{eqnarray*}
which ends the proof.
\end{Proof}

\footnotesize

\def\cprime{$'$} \def\polhk#1{\setbox0=\hbox{#1}{\ooalign{\hidewidth
  \lower1.5ex\hbox{`}\hidewidth\crcr\unhbox0}}}
  \def\polhk#1{\setbox0=\hbox{#1}{\ooalign{\hidewidth
  \lower1.5ex\hbox{`}\hidewidth\crcr\unhbox0}}} \def\cprime{$'$}

\end{document}